# Can Decentralized Algorithms Outperform Centralized Algorithms? A Case Study for Decentralized Parallel Stochastic Gradient Descent

Xiangru Lian[*1], Ce Zhang[†2], Huan Zhang[‡3], Cho-Jui Hsieh[§3], Wei Zhang[¶4], and Ji Liu[∥1,5]

[1]University of Rochester
[2]ETH Zurich
[3]University of California, Davis
[4]IBM T. J. Watson Research Center
[5]Tencent AI lab

## Abstract

*Most distributed machine learning systems nowadays, including TensorFlow and CNTK, are built in a centralized fashion. One bottleneck of centralized algorithms lies on high communication cost on the central node. Motivated by this, we ask, can decentralized algorithms be faster than its centralized counterpart?*

*Although decentralized PSGD (D-PSGD) algorithms have been studied by the control community, existing analysis and theory do not show any advantage over centralized PSGD (C-PSGD) algorithms, simply assuming the application scenario where only the decentralized network is available. In this paper, we study a D-PSGD algorithm and provide the first theoretical analysis that indicates a regime in which decentralized algorithms might outperform centralized algorithms for distributed stochastic gradient descent. This is because D-PSGD has comparable total computational complexities to C-PSGD but requires much less communication cost on the busiest node. We further conduct an empirical study to validate our theoretical analysis across multiple frameworks (CNTK and Torch), different network configurations, and computation platforms up to 112 GPUs. On network configurations with low bandwidth or high latency, D-PSGD can be up to one order of magnitude faster than its well-optimized centralized counterparts.*

## 1 Introduction

In the context of distributed machine learning, decentralized algorithms have long been treated as a *compromise* — when the underlying network topology does not allow centralized communication, one *has*

[*]xiangru@yandex.com
[†]ce.zhang@inf.ethz.ch
[‡]victzhang@gmail.com
[§]chohsieh@ucdavis.edu
[¶]weiz@us.ibm.com
[∥]ji.liu.uwisc@gmail.com

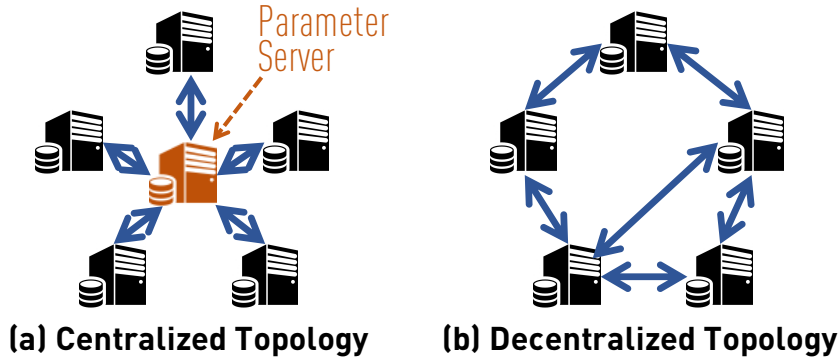

Figure 1: An illustration of different network topologies.

*to* resort to decentralized communication, while, understandably, pay for the "cost of being decentralized". In fact, most distributed machine learning systems nowadays, including TensorFlow and CNTK, are built in a centralized fashion. But *can decentralized algorithms be faster than its centralized counterpart?* In this paper, we provide the first theoretical analysis, verified by empirical experiments, for a positive answer to this question.

We consider solving the following stochastic optimization problem

$$\min_{x \in \mathbb{R}^N} f(x) := \mathbb{E}_{\xi \sim \mathcal{D}} F(x; \xi), \tag{1}$$

where $\mathcal{D}$ is a predefined distribution and $\xi$ is a random variable usually referring to a data sample in machine learning. This formulation summarizes many popular machine learning models including deep learning [LeCun et al., 2015], linear regression, and logistic regression.

Parallel stochastic gradient descent (PSGD) methods are leading algorithms in solving large-scale machine learning problems such as deep learning [Dean et al., 2012, Li et al., 2014], matrix completion [Recht et al., 2011, Zhuang et al., 2013] and SVM. Existing PSGD algorithms are mostly designed for centralized network topology, for example, parameter server [Li et al., 2014], where there is a central node connected with multiple nodes as shown in Figure 1(a). The central node aggregates the stochastic gradients computed from all other nodes and updates the model parameter, for example, the weights of a neural network. The potential bottleneck of the centralized network topology lies on the communication traffic jam on the central node, because all nodes need to communicate with it concurrently iteratively. The performance will be significantly degraded when the network bandwidth is low.[1] These motivate us to study algorithms for *decentralized* topologies, where all nodes can only communicate with its neighbors and there is no such a central node, shown in Figure 1(b).

Although decentralized algorithms have been studied as consensus optimization in the control community and used for preserving data privacy [Ram et al., 2009a, Yan et al., 2013, Yuan et al., 2016], for the application scenario where only the decentralized network is available, it is still an open question **if decentralized methods could have advantages over centralized algorithms** in some scenarios in case both types of communication patterns are feasible — for example, on a supercomputer with thousands of nodes, *should we use decentralized or centralized communication?* Existing theory and analysis either do not make such comparison [Bianchi et al., 2013, Ram et al., 2009a, Srivastava and Nedic, 2011, Sundhar Ram et al., 2010] or implicitly indicate that decentralized algorithms were much worse than centralized algorithms in terms of computational complexity and total communication complexity [Aybat et al., 2015, Lan et al., 2017, Ram et al., 2010, Zhang and Kwok, 2014]. This paper gives a positive result for decentralized algorithms by studying a decentralized PSGD (D-PSGD) algorithm on the connected decentralized network. Our theory indicates that D-PSGD admits similar total computational complexity but requires much less communication for the busiest node. Table 1 shows a quick comparison between C-PSGD and D-PSGD with respect to the computation and communication complexity. Our contributions are:

- We theoretically justify the potential advantage of decentralizedalgorithms over centralized algorithms. Instead of treating decentralized algorithms as a compromise one has to make, we are the first to conduct a theoretical analysis that identifies cases in which decentralized algorithms can be faster than its centralized counterpart.

- We theoretically analyze the scalability behavior of decentralized SGD when more nodes are used. Surprisingly, we show that, when more nodes are available, decentralized algorithms can bring speedup, asymptotically linearly, with respect to computational complexity. To our best knowledge, this is the first speedup result related to decentralized algorithms.

- We conduct an empirical study to validate our theoretical analysis of D-PSGD and different C-PSGD variants (e.g., plain SGD, EASGD [Zhang et al., 2015]). We observe similar computational complexity as our theory indicates; on networks with low bandwidth or high latency, D-PSGD can be up to $10\times$ faster than C-PSGD. Our result holds across multiple frameworks (CNTK and Torch), different network configurations, and computation platforms up to 112 GPUs. This indicates promising future direction in pushing the research horizon of machine learning systems from pure centralized topology to a more decentralized fashion.

---

[1]There has been research in how to accommodate this problem by having multiple parameter servers communicating with efficient MPI ALLREDUCE primitives. As we will see in the experiments, these methods, on the other hand, might suffer when the network latency is high.

| Algorithm | communication complexity on the busiest node | computational complexity |
|---|---|---|
| C-PSGD (mini-batch SGD) | $O(n)$ | $O(\frac{n}{\epsilon} + \frac{1}{\epsilon^2})$ |
| D-PSGD | $O(\text{Deg(network)})$ | $O(\frac{n}{\epsilon} + \frac{1}{\epsilon^2})$ |

Table 1: Comparison of C-PSGD and D-PSGD. The unit of the communication cost is the number of stochastic gradients or optimization variables. $n$ is the number of nodes. The computational complexity is the number of stochastic gradient evaluations we need to get a $\epsilon$-approximation solution, which is defined in (3).

**Definitions and notations** Throughout this paper, we use following notation and definitions:

- $\|\cdot\|$ denotes the vector $\ell_2$ norm or the matrix spectral norm depending on the argument.
- $\|\cdot\|_F$ denotes the matrix Frobenius norm.
- $\nabla f(\cdot)$ denotes the gradient of a function $f$.
- $\mathbf{1}_n$ denotes the column vector in $\mathbb{R}^n$ with 1 for all elements.
- $f^*$ denotes the optimal solution of (1).
- $\lambda_i(\cdot)$ denotes the $i$-th largest eigenvalue of a matrix.

## 2  Related work

In the following, we use *K* and *n* to refer to the number of iterations and the number of nodes.

**Stochastic Gradient Descent (SGD)**  SGD is a powerful approach for solving large scale machine learning. The well known convergence rate of stochastic gradient is $O(1/\sqrt{K})$ for convex problems and $O(1/K)$ for strongly convex problems [Moulines and Bach, 2011, Nemirovski et al., 2009]. SGD is closely related to online learning algorithms, for example, Crammer et al. [2006], Shalev-Shwartz [2011], Yang et al. [2014]. For SGD on nonconvex optimization, an ergodic convergence rate of $O(1/\sqrt{K})$ is proved in Ghadimi and Lan [2013].

**Centralized parallel SGD**  For CENTRALIZED PARALLEL SGD (C-PSGD) algorithms, the most popular implementation is the parameter server implementation, which is essentially mini-batch SGD admitting a convergence rate of $O(1/\sqrt{Kn})$ [Agarwal and Duchi, 2011, Dekel et al., 2012, Lian et al., 2015], where in each iteration *n* stochastic gradients are evaluated. In this implementation there is a parameter server communicating with all nodes. The linear speedup is implied by the convergence rate automatically. More implementation details for C-PSGD can be found in Chen et al. [2016], Dean et al. [2012], Li et al. [2014], Zinkevich et al. [2010]. The asynchronous version of centralized parallel SGD is proved to maintain a linear speedup on convex, strongly convex and nonconvex objectives when the staleness of the gradient is bounded [Agarwal and Duchi, 2011, Feyzmahdavian et al., 2015, Lian et al., 2015, 2016, Recht et al., 2011].

**Decentralized parallel stochastic algorithms**  There are existing work on decentralized parallel stochastic gradient where there is no central node (parameter server). They look similar to D-PSGD but *none of them is proved to have speedup when we increase the number of nodes*. For example, Lan et al. [2017] gave a decentralized stochastic algorithm with a computational complexity of $O(n/\epsilon^2)$ for general convex objectives and $O(n/\epsilon)$ for strongly convex objectives. Sirb and Ye [2016] gave a $O(n/\epsilon^2)$ complexity for convex objectives with an asynchronous decentralized stochastic algorithm. These bounds for the complexity are proportional to *n*, which means no speedup is shown. We review other related work in the following.

An algorithm similar to D-PSGD in both synchronous and asynchronous fashion was studied in Ram et al. [2009a, 2010], Srivastava and Nedic [2011], Sundhar Ram et al. [2010]. The difference is that in their algorithm a node cannot do communication and computation simultaneously. The algorithm in Srivastava and Nedic [2011] optimizes the convex objective, however, to derive an error bound it requires bounded domain and each term in the objective function to be strongly convex. Sundhar Ram et al. [2010] is its subgradient variant. The analysis in Ram et al. [2009a, 2010], Srivastava and Nedic [2011], Sundhar Ram et al. [2010] requires the gradients of each term of the objective to be bounded by a constant. The analysis in Bianchi et al. [2013] uses strong non-standard assumptions for a decentralized stochastic algorithm, which requires continuously increasing communication cost when we run the algorithm since the second largest eigenvalue of the averaging matrix needs to be decreasing to 0 when the algorithm is running.

**Other decentralized algorithms**  In other areas including control, privacy and wireless sensing network, there are work on the consensus problem for which decentralized algorithms are studied to compute the mean of all the data distributed on multiple nodes [Aysal et al., 2009, Boyd et al., 2005, Carli et al., 2010, Fagnani and Zampieri, 2008, Olfati-Saber et al., 2007, Schenato and Gamba, 2007]. Lu et al. [2010] showed a gossip algorithm applied on convex objectives converges to the solution but no convergence rate was

**Algorithm 1** Decentralized Parallel Stochastic Gradient Descent (D-PSGD) on the $i$th node
___
**Require:** initial point $x_{0,i} = x_0$, step length $\gamma$, weight matrix $W$, and number of iterations $K$
1: **for** $k = 0, 1, 2, \ldots, K - 1$ **do**
2:     Randomly sample $\xi_{k,i}$ from local data of the $i$-th node
3:     Compute a local stochastic gradient based on $\xi_{k,i}$ and current optimization variable $x_{k,i}$: $\nabla F_i(x_{k,i}; \xi_{k,i})$ [a]
4:     Compute the neighborhood weighted average by fetching optimization variables from neighbors: $x_{k+\frac{1}{2},i} = \sum_{j=1}^n W_{ij} x_{k,j}$ [b]
5:     Update the local optimization variable $x_{k+1,i} \leftarrow x_{k+\frac{1}{2},i} - \gamma \nabla F_i(x_{k,i}; \xi_{k,i})$ [c]
6: **end for**
7: **Output:** $\frac{1}{n} \sum_{i=1}^n x_{K,i}$ [d]

___
[a] Note that the stochastic gradient computed in can be replaced with a mini-batch of stochastic gradients, which will not hurt our theoretical results.

[b] Note that the Line 3 and Line 4 can be run in parallel.

[c] Note that the Line 4 and step Line 5 can be exchanged. That is, we first update the local stochastic gradient into the local optimization variable, and then average the local optimization variable with neighbors. This does not hurt our theoretical analysis. When Line 4 is logically before Line 5, then Line 3 and Line 4 can be run in parallel. That is to say, if the communication time used by Line 4 is smaller than the computation time used by Line 3, the communication time can be completely hidden (it is overlapped by the computation time).

[d] We will prove that the local optimization variables in the nodes will converge together, so it is also safe to use the local optimization variable of a single node as an estimation of the solution.

___

shown. Mokhtari and Ribeiro [2016] analyzed decentralized SAG and SAGA algorithms. They are not shown to have speedup, and a table of all stochastic gradients need to be saved in the storage or memory. Decentralized gradient descent on convex and strongly convex problems was analyzed in Yuan et al. [2016]. Nedic and Ozdaglar [2009], Ram et al. [2009b] are similar to Yuan et al. [2016] but they use subgradients. The algorithm in Nedic and Ozdaglar [2009], Ram et al. [2009b], Yuan et al. [2016] does not converge to the exact solution due to the inconsistent nature of decentralized gradient descent. This was fixed by Shi et al. [2015] using a modified algorithm. Wu et al. [2016] analyzed an asynchronous version of decentralized gradient descent with some modification like in Shi et al. [2015] and showed the algorithm converges to a solution when $K \to \infty$. Aybat et al. [2015], Shi et al., Zhang and Kwok [2014] analyzed decentralized ADMM algorithms and they are not shown to have speedup. From all of these reviewed papers, it is still unclear if decentralized algorithms can outperform centralized algorithms.

## 3 Decentralized parallel stochastic gradient descent (D-PSGD)

This section introduces the D-PSGD algorithm. We represent the decentralized communication topology with an undirected graph with weights: $(V, W)$. $V$ denotes the set of $n$ computational nodes: $V := \{1, 2, \cdots, n\}$. $W \in \mathbb{R}^{n \times n}$ is a symmetric doubly stochastic matrix, which means (i) $W_{ij} \in [0,1], \forall i, j$, (ii) $W_{ij} = W_{ji}$ for all $i, j$, and (ii) $\sum_j W_{ij} = 1$ for all $i$. We use $W_{ij}$ to encode how much node $j$ can affect node $i$, while $W_{ij} = 0$ means node $i$ and $j$ are disconnected.

To design distributed algorithms on a decentralized network, we first distribute the data onto all nodes such that the original objective defined in (1) can be rewritten into

$$\min_{x \in \mathbb{R}^N} \quad f(x) = \frac{1}{n} \sum_{i=1}^n \underbrace{\mathbb{E}_{\xi \sim \mathcal{D}_i} F_i(x; \xi)}_{=: f_i(x)}. \tag{2}$$

There are two simple ways to achieve (2), both of which can be captured by our theoretical analysis and they both imply $F_i(\cdot;\cdot) = F(\cdot;\cdot), \forall i$.

**Strategy-1** All distributions $\mathcal{D}_i$'s are the same as $\mathcal{D}$, that is, all nodes can access a shared database;

**Strategy-2** $n$ nodes partition all data in the database and appropriately define a distribution for sampling local data, for example, if $\mathcal{D}$ is the uniform distribution over all data, $\mathcal{D}_i$ can be defined to be the uniform distribution over local data.

The D-PSGD algorithm is a synchronous parallel algorithm. All nodes are usually synchronized by a clock. Each node maintains its own local variable and runs the protocol in Algorithm 1 concurrently, which includes three key steps at iterate $k$:

- Each node computes the stochastic gradient $\nabla F_i(x_{k,i};\xi_{k,i})$[2] using the current local variable $x_{k,i}$, where $k$ is the iterate number and $i$ is the node index;

- When the synchronization barrier is met, each node exchanges local variables with its neighbors and average the local variables it receives with its own local variable;

- Each node update its local variable using the average and the local stochastic gradient.

To view the D-PSGD algorithm from a global view, at iterate $k$, we define the concatenation of all local variables, random samples, stochastic gradients by matrix $X_k \in \mathbb{R}^{N\times n}$, vector $\xi_k \in \mathbb{R}^n$, and $\partial F(X_k, \xi_k)$, respectively:

$$X_k := \begin{bmatrix} x_{k,1} & \cdots & x_{k,n} \end{bmatrix} \in \mathbb{R}^{N\times n}, \quad \xi_k := \begin{bmatrix} \xi_{k,1} & \cdots & \xi_{k,n} \end{bmatrix}^\top \in \mathbb{R}^n,$$

$$\partial F(X_k, \xi_k) := \begin{bmatrix} \nabla F_1(x_{k,1};\xi_{k,1}) & \nabla F_2(x_{k,2};\xi_{k,2}) & \cdots & \nabla F_n(x_{k,n};\xi_{k,n}) \end{bmatrix} \in \mathbb{R}^{N\times n}.$$

Then the $k$-th iterate of Algorithm 1 can be viewed as the following update

$$X_{k+1} \leftarrow X_k W - \gamma \partial F(X_k; \xi_k).$$

We say the algorithm gives an $\epsilon$-approximation solution if

$$K^{-1}\left(\sum_{k=0}^{K-1} \mathbb{E}\left\|\nabla f\left(\frac{X_k \mathbf{1}_n}{n}\right)\right\|^2\right) \leq \epsilon. \tag{3}$$

## 4 Convergence rate analysis

This section provides the analysis for the convergence rate of the D-PSGD algorithm. Our analysis will show that the convergence rate of D-PSGD w.r.t. iterations is similar to the C-PSGD (or mini-batch SGD) [Agarwal and Duchi, 2011, Dekel et al., 2012, Lian et al., 2015], but D-PSGD avoids the communication traffic jam on the parameter server.

To show the convergence results, we first define

$$\partial f(X_k) := \begin{bmatrix} \nabla f_1(x_{k,1}) & \nabla f_2(x_{k,2}) & \cdots & \nabla f_n(x_{k,n}) \end{bmatrix} \in \mathbb{R}^{N\times n},$$

where functions $f_i(\cdot)$'s are defined in (2).

---
[2] It can be easily extended to mini-batch stochastic gradient descent.

**Assumption 1.** *Throughout this paper, we make the following commonly used assumptions:*

1. **Lipschitzian gradient:** *All function $f_i(\cdot)$'s are with L-Lipschitzian gradients.*

2. **Spectral gap:** *Given the symmetric doubly stochastic matrix W, we define $\rho := (\max\{|\lambda_2(W)|, |\lambda_n(W)|\})^2$. We assume $\rho < 1$.*

3. **Bounded variance:** *Assume the variance of stochastic gradient*
$$\mathbb{E}_{i \sim \mathcal{U}([n])} \mathbb{E}_{\xi \sim \mathcal{D}_i} \|\nabla F_i(x; \xi) - \nabla f(x)\|^2$$
*is bounded for any x with i uniformly sampled from $\{1, \ldots, n\}$ and $\xi$ from the distribution $\mathcal{D}_i$. This implies there exist constants $\sigma, \varsigma$ such that*
$$\mathbb{E}_{\xi \sim \mathcal{D}_i} \|\nabla F_i(x; \xi) - \nabla f_i(x)\|^2 \leq \sigma^2, \forall i, \forall x,$$
$$\mathbb{E}_{i \sim \mathcal{U}([n])} \|\nabla f_i(x) - \nabla f(x)\|^2 \leq \varsigma^2, \forall x.$$
*Note that if all nodes can access the shared database, then $\varsigma = 0$.*

4. **Start from 0:** *We assume $X_0 = 0$. This assumption simplifies the proof w.l.o.g.*

Let
$$D_1 := \left(\frac{1}{2} - \frac{9\gamma^2 L^2 n}{(1-\sqrt{\rho})^2 D_2}\right), \quad D_2 := \left(1 - \frac{18\gamma^2}{(1-\sqrt{\rho})^2} nL^2\right).$$
Under Assumption 1, we have the following convergence result for Algorithm 1.

**Theorem 1** (Convergence of Algorithm 1). *Under Assumption 1, we have the following convergence rate for Algorithm 1:*

$$\frac{1}{K}\left(\frac{1-\gamma L}{2} \sum_{k=0}^{K-1} \mathbb{E} \left\|\frac{\partial f(X_k)\mathbf{1}_n}{n}\right\|^2 + D_1 \sum_{k=0}^{K-1} \mathbb{E} \left\|\nabla f\left(\frac{X_k \mathbf{1}_n}{n}\right)\right\|^2\right)$$
$$\leq \frac{f(0) - f^*}{\gamma K} + \frac{\gamma L}{2n}\sigma^2 + \frac{\gamma^2 L^2 n \sigma^2}{(1-\rho) D_2} + \frac{9\gamma^2 L^2 n \varsigma^2}{(1-\sqrt{\rho})^2 D_2}.$$

Noting that $\frac{X_k \mathbf{1}_n}{n} = \frac{1}{n}\sum_{i=1}^n x_{k,i}$, this theorem characterizes the convergence of the average of all local optimization variables $x_{k,i}$. To take a closer look at this result, we appropriately choose the step length in Theorem 1 to obtain the following result:

**Corollary 2.** *Under the same assumptions as in Theorem 1, if we set $\gamma = \frac{1}{2L + \sigma\sqrt{K/n}}$, for Algorithm 1 we have the following convergence rate:*

$$\frac{\sum_{k=0}^{K-1} \mathbb{E} \left\|\nabla f\left(\frac{X_k \mathbf{1}_n}{n}\right)\right\|^2}{K} \leq \frac{8(f(0) - f^*)L}{K} + \frac{(8f(0) - 8f^* + 4L)\sigma}{\sqrt{Kn}}. \quad (4)$$

*if the total number of iterate K is sufficiently large, in particular,*

$$K \geq \frac{4L^4 n^5}{\sigma^6 (f(0) - f^* + L)^2}\left(\frac{\sigma^2}{1-\rho} + \frac{9\varsigma^2}{(1-\sqrt{\rho})^2}\right)^2, \text{ and} \quad (5)$$

$$K \geq \frac{72 L^2 n^2}{\sigma^2 (1-\sqrt{\rho})^2}. \quad (6)$$

This result basically suggests that the convergence rate for D-PSGD is $O\left(\frac{1}{K} + \frac{1}{\sqrt{nK}}\right)$, if $K$ is large enough. We highlight two key observations from this result:

**Linear speedup** When $K$ is large enough, the $\frac{1}{K}$ term will be dominated by the $\frac{1}{\sqrt{Kn}}$ term which leads to a $\frac{1}{\sqrt{nK}}$ convergence rate. It indicates that the total computational complexity [3] to achieve an $\epsilon$-approximation solution (3) is bounded by $O\left(\frac{1}{\epsilon^2}\right)$. Since the total number of nodes does not affect the total complexity, a single node only shares a computational complexity of $O\left(\frac{1}{n\epsilon^2}\right)$. Thus linear speedup can be achieved by D-PSGD asymptotically w.r.t. computational complexity.

**D-PSGD can be better than C-PSGD** Note that this rate is the same as C-PSGD (or mini-batch SGD with mini-batch size $n$) [Agarwal and Duchi, 2011, Dekel et al., 2012, Lian et al., 2015]. The advantage of D-PSGD over C-PSGD is to avoid the communication traffic jam. At each iteration, the *maximal* communication cost for every single node is $O$(the degree of the network) for D-PSGD, in contrast with $O(n)$ for C-PSGD. The degree of the network could be much smaller than $O(n)$, e.g., it could be $O(1)$ in the special case of a ring.

The key difference from most existing analysis for decentralized algorithms lies on that we do not use the boundedness assumption for domain or gradient or stochastic gradient. Those boundedness assumptions can significantly simplify the proof but lose some subtle structures in the problem.

The linear speedup indicated by Corollary 4 requires the total number of iteration $K$ is sufficiently large. The following special example gives a concrete bound of $K$ for the ring network topology.

**Theorem 3.** (**Ring network**) *Choose the steplength $\gamma$ in the same as Corollary 2 and consider the ring network topology with corresponding W in the form of*

$$W = \begin{pmatrix} 1/3 & 1/3 & & & & & 1/3 \\ 1/3 & 1/3 & 1/3 & & & & \\ & 1/3 & 1/3 & \ddots & & & \\ & & \ddots & \ddots & 1/3 & & \\ & & & 1/3 & 1/3 & 1/3 & \\ 1/3 & & & & & 1/3 & 1/3 \end{pmatrix} \in \mathbb{R}^{n \times n}.$$

*Under Assumption 1, Algorithm 1 achieves the same convergence rate in* (4), *which indicates a linear speedup can be achieved, if the number of involved nodes is bounded by*

- $n = O(K^{1/9})$, *if apply **strategy-1** distributing data ($\varsigma = 0$);*
- $n = O(K^{1/13})$, *if apply **strategy-2** distributing data ($\varsigma > 0$),*

*where the capital "O" swallows $\sigma, \varsigma, L$, and $f(0) - f^*$.*

This result considers a special decentralized network topology: ring network, where each node can only exchange information with its two neighbors. The linear speedup can be achieved up to $K^{1/9}$ and $K^{1/13}$ for different scenarios. These two upper bound can be improved potentially. This is the first work to show the speedup for decentralized algorithms, to the best of our knowledge.

---
[3]The computation of a single stochastic gradient counts 1. So the computational complexity per iteration is $O(n)$.

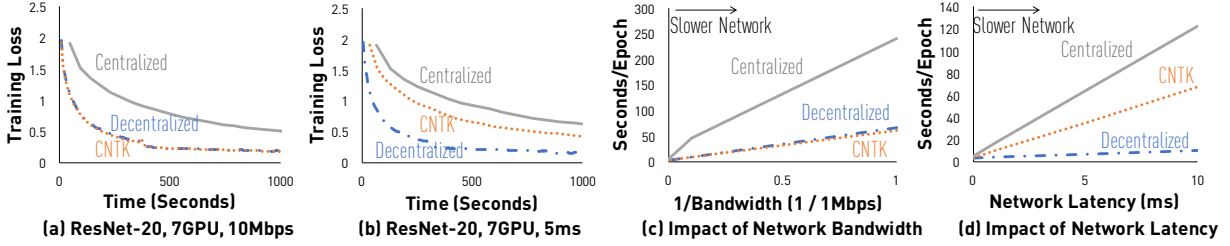

Figure 2: Comparison between D-PSGD and two centralized implementations (7 and 10 GPUs).

In this section, we mainly investigate the convergence rate for the average of all local variables $\{x_{k,i}\}_{i=1}^{n}$. Actually one can also obtain a similar rate for each individual $x_{k,i}$, since all nodes achieve the consensus quickly, in particular, the running average of $\mathbb{E}\left\|\frac{\sum_{i'=1}^{n} x_{k,i'}}{n} - x_{k,i}\right\|^2$ converges to 0 with a $O(1/K)$ rate, where the "$O$" swallows $n, \rho, \sigma, \varsigma, L$ and $f(0) - f^*$. This result can be formally summarized into the following theorem:

**Theorem 4.** *With $\gamma = \frac{1}{2L + \sigma\sqrt{K/n}}$ under the same assumptions as in Corollary 2 we have*

$$(Kn)^{-1} \mathbb{E} \sum_{k=0}^{K-1} \sum_{i=1}^{n} \left\|\frac{\sum_{i'=1}^{n} x_{k,i'}}{n} - x_{k,i}\right\|^2 \leqslant n\gamma^2 \frac{A}{D_2},$$

*where*

$$A := \frac{2\sigma^2}{1-\rho} + \frac{18\varsigma^2}{(1-\sqrt{\rho})^2} + \frac{L^2}{D_1}\left(\frac{\sigma^2}{1-\rho} + \frac{9\varsigma^2}{(1-\sqrt{\rho})^2}\right) \\ + \frac{18}{(1-\sqrt{\rho})^2}\left(\frac{f(0) - f^*}{\gamma K} + \frac{\gamma L \sigma^2}{2nD_1}\right).$$

Choosing $\gamma$ in the way in Corollary 4, we can see that the consensus will be achieved in the rate $O(1/K)$.

## 5 Experiments

We validate our theory with experiments that compared D-PSGD with other centralized implementations. We run experiments on clusters up to 112 GPUs and show that, on some network configurations, D-PSGD can outperform well-optimized centralized implementations by an order of magnitude.

### 5.1 Experiment setting

**Datasets and models** We evaluate D-PSGD on two machine learning tasks, namely (1) image classification, and (2) Natural Language Processing (NLP). For image classification we train ResNet [He et al., 2015] with different number of layers on CIFAR-10 [Krizhevsky, 2009]; for speech recognition, we train both proprietary and public dataset on a proprietary CNN model that we get from our industry partner [Feng et al., 2016, Lin et al., 2017]. We leave the result of NLP to the supplementary material because the results are similar to that of image classification.

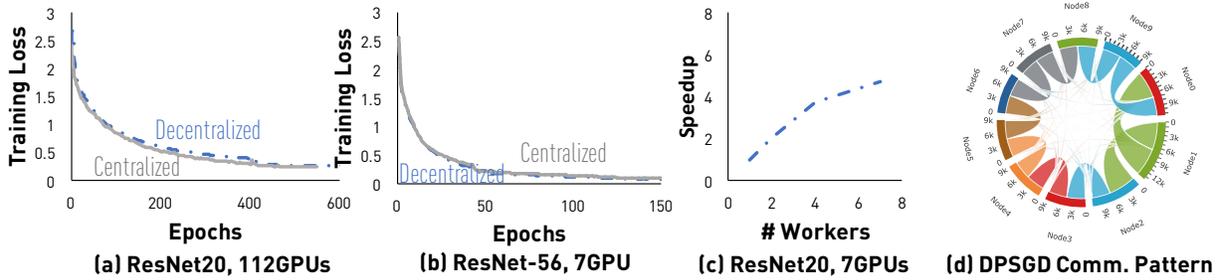

Figure 3: (a) Convergence Rate; (b) D-PSGD Speedup; (c) D-PSGD Communication Patterns.

**Implementations and setups** We implement D-PSGD on two different frameworks, namely Microsoft CNTK and Torch. We evaluate four SGD implementations:

1. **CNTK.** We compare with the standard CNTK implementation of synchronous SGD. The implementation is based on MPI's AllReduce primitive.

2. **Centralized.** We implemented the standard parameter server-based synchronous SGD using MPI. One node will serve as the parameter server in our implementation.

3. **Decentralized.** We implemented our D-PSGD algorithm using MPI within CNTK.

4. **EASGD.** We compare with the standard EASGD implementation of Torch.

All three implementations are compiled with gcc 7.1, cuDNN 5.0, OpenMPI 2.1.1. We fork from CNTK after commit `57d7b9d` and enable distributed minibatch reading for all of our experiments.

During training, we keep the local batch size of each node the same as the reference configurations provided by CNTK. We tune learning rate for each SGD variant and report the best configuration.

**Machines/Clusters** We conduct experiments on three different machines/clusters:

1. **7GPUs.** A single local machine with 8 GPUs, each of which is a Nvidia TITAN Xp.

2. **10GPUs.** 10 `p2.xlarge` EC2 instances, each of which has one Nvidia K80 GPU.

3. **16GPUs.** 16 local machines, each of which has two Xeon E5-2680 8-core processors and a NVIDIA K20 GPU. Machines are connected by Gigabit Ethernet in this case.

4. **112GPUs.** 4 `p2.16xlarge` and 6 `p2.8xlarge` EC2 instances. Each `p2.16xlarge` (resp. `p2.8xlarge`) instance has 16 (resp. 8) Nvidia K80 GPUs.

In all of our experiments, we use each GPU as a node.

## 5.2 Results on CNTK

**End-to-end performance**  We first validate that, under certain network configurations, D-PSGD converges faster, in wall-clock time, to a solution that has the same quality of centralized SGD. Figure 2(a, b) and Figure 3(a) shows the result of training ResNet20 on 7GPUs. We see that D-PSGD converges faster than both centralized SGD competitors. This is because when the network is slow, both centralized SGD competitors take more time per epoch due to communication overheads. Figure 3(a, b) illustrates the convergence with respect to the number of epochs, and D-PSGD shows similar convergence rate as centralized SGD even with 112 nodes.

**Speedup**  The end-to-end speedup of D-PSGD over centralized SGD highly depends on the underlying network. We use the `tc` command to manually vary the network bandwidth and latency and compare the wall-clock time that all three SGD implementations need to finish one epoch.

Figure 2(c, d) shows the result. We see that, when the network has high bandwidth and low latency, not surprisingly, all three SGD implementations have similar speed. This is because in this case, the communication is never the system bottleneck. However, when the bandwidth becomes smaller (Figure 2(c)) or the latency becomes higher (Figure 2(d)), both centralized SGD implementations slow down significantly. In some cases, D-PSGD can be even one order of magnitude faster than its centralized competitors. Compared with **Centralized** (implemented with a parameter server), D-PSGD has more balanced communication patterns between nodes and thus outperforms **Centralized** in low-bandwidth networks; compared with **CNTK** (implemented with AllReduce), D-PSGD needs fewer number of communications between nodes and thus outperforms **CNTK** in high-latency networks. Figure 3(c) illustrates the communication between nodes for one run of D-PSGD.

We also vary the number of GPUs that D-PSGD uses and report the speed up over a single GPU to reach the same loss. Figure 3(b) shows the result on a machine with 7GPUs. We see that, up to 4 GPUs, D-PSGD shows near linear speed up. When all seven GPUs are used, D-PSGD achieves up to $5\times$ speed up. This subliner speed up for 7 GPUs is due to the synchronization cost but also that our machine only has 4 PCIe channels and thus more than two GPUs will share PCIe bandwidths.

## 5.3 Results on Torch

We provide report results for the experiment of D-PSGD and EASGD. For this set of experiments we use a 32-layer residual network and CIFAR-10 dataset. We use up to 16 machines, and each machine includes two Xeon E5-2680 8-core processors and a NVIDIA K20 GPU. Worker machines are connected in a logical ring as described in Theorem 3. Connections between D-PSGD nodes are made via TCP socks, and EASGD uses MPI for communication. Because D-PSGD do not have a centralized model, we average all models from different machines as our final model to evaluate. In practical training, this only needs to be done after the last epoch with an all-reduce operation. For EASGD, we evaluate the central model on the parameter server.

One remarkable feature of this experiment is that we use inexpensive Gigabit Ethernet to connect all machines, and we are able to practically observe network congestion with centralized parameter server approach, even with a relatively small (ResNet-32) model. Although in practice, network with much higher bandwidth are available (e.g., InfiniBand), we also want to use larger model or more machines, so that network bandwidth can always become a bottleneck. We practically show that D-PSGD has better

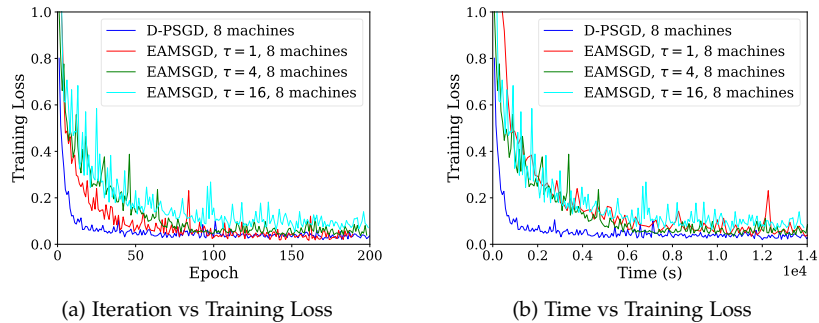

(a) Iteration vs Training Loss  (b) Time vs Training Loss

Figure 4: Convergence comparison between D-PSGD and EAMSGD (EASGD's momentum variant).

scalability than centralized approaches when network bandwidth becomes a constraint.

**Comparison to EASGD**   Elastic Averaging SGD (EASGD) [Zhang et al., 2015] is an improved parameter server approach that outperforms traditional parameter server [Dean et al., 2012]. It makes each node perform more exploration by allowing local parameters to fluctuate around the central variable. We add ResNet-32 [He et al., 2016] with CIFAR-10 into the EASGD's Torch experiment code[4] and also implement our algorithm in Torch. Both algorithms run at the same speed on a single GPU so there is no implementation bias. Unlike the previous experiment which uses high bandwidth PCI-e or 10Gbits network for inter-GPU communication, we use 9 physical machines (1 as parameter server) with a single K20 GPU each, connected by inexpensive Gigabit Ethernet. For D-PSGD we use a logical ring connection between nodes as in Theorem 3. For EASGD we set moving rate $\beta = 0.9$ and use its momentum variant (EAMSGD). For both algorithms we set learning rate to 0.1, momentum to 0.9. $\tau = \{1, 4, 16\}$ is a hyper-parameter in EASGD controlling the number of mini-batches before communicating with the server.

Figure 4 shows that D-PSGD outperforms EASGD with a large margin in this setting. EASGD with $\tau = 1$ has good convergence, but its large bandwidth requirement saturates the network and slows down nodes. When $\tau = 4, 16$ EASGD converges slower than D-PSGD as there is less communication. D-PSGD allows more communication in an efficient way without reaching the network bottleneck. Moreover, D-PSGD is synchronous and shows less convergence fluctuation comparing with EASGD.

**Accuracy comparison with EASGD**   We have shown the training loss comparison between D-PSGD and EASGD, and we now show additional figures comparing training error and test error in our experiment, as in Figure 5 and 6. We observe similar results as we have seen in section 5.3; D-PSGD can achieve good accuracy noticeably faster than EASGD.

**Scalability of D-PSGD**   In this experiment, we run D-PSGD on 1, 4, 8, 16 machines and compare convergence speed and error. For experiments involving 16 machines, each machine also connects to one additional machine which has the largest topological distance on the ring besides its two logical neighbours. We found that this can help information flow and get better convergence.

In Figure 10, 11 and 12 we can observe that D-PSGD scales very well when the number of machines is growing. Also, comparing with the single machine SGD, D-PSGD has minimum overhead; we measure the per-epoch training time only increases by 3% comparing to single machine SGD, but D-PSGD's convergence speed is much faster. To reach a training loss of 0.2, we need about 80 epochs with 1 machine, 20 epochs with 4 machines, 10 epochs with 8 machines and only 5 epochs with 16 machines. The observed linear speedup justifies the correctness of our theory.

---
[4]https://github.com/sixin-zh/mpiT.git

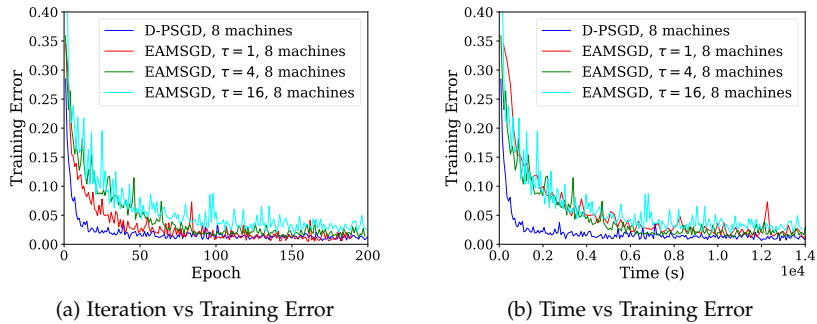

(a) Iteration vs Training Error  (b) Time vs Training Error

Figure 5: Training Error comparison between D-PSGD and EAMSGD (EASGD's momentum variant)

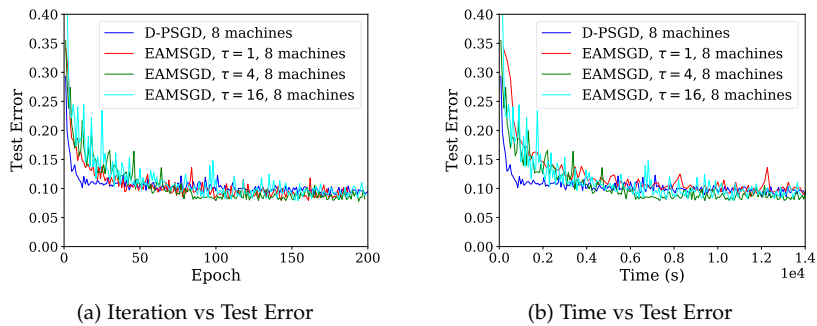

(a) Iteration vs Test Error  (b) Time vs Test Error

Figure 6: Test Error comparison between D-PSGD and EAMSGD (EASGD's momentum variant)

**Generalization ability of D-PSGD**  In our previous experiments we set the learning rate to fixed 0.1. To complete Residual network training, we need to decrease the learning rate after some epochs. We follow the learning rate schedule in ResNet paper [He et al., 2016], and decrease the learning rate to 0.01 at epoch 80. We observe training/test loss and error, as shown in figure 10, 11 and 12. For D-PSGD, we can tune a better learning rate schedule, but parameter tuning is not the focus of our experiments; rather, we would like to see if D-PSGD can achieve the same best ResNet accuracy as reported by the literature.

The test error of D-PSGD after 160 epoch is 0.0715, 0.0746 and 0.0735, for 4, 8 and 16 machines, respectively. He et al. [2016] reports 0.0751 error for the same 32-layer residual network, and we can reliably outperform the reported error level regardless of different numbers of machines used. Thus, D-PSGD does not negatively affect (or perhaps helps) generalization.

**Network utilization**  During the experiment, we measure the network bandwidth on each machine. Because every machine is identical on the network, the measured bandwidth are the same on each machines. For experiment with 4 and 8 machines, the required bandwidth is about 22 MB/s. With 16 machines the required bandwidth is about 33 MB/s because we have an additional link. The required bandwidth is related to GPU performance; if GPU can compute each minibatch faster, the required bandwidth also increases proportionally. Considering the practical bandwidth of Gigabit Ethernet is about 100 ~120 MB/s, Our algorithm can handle a 4 ~5 times faster GPU (or GPUs) easily, even with an inexpensive gigabit connection.

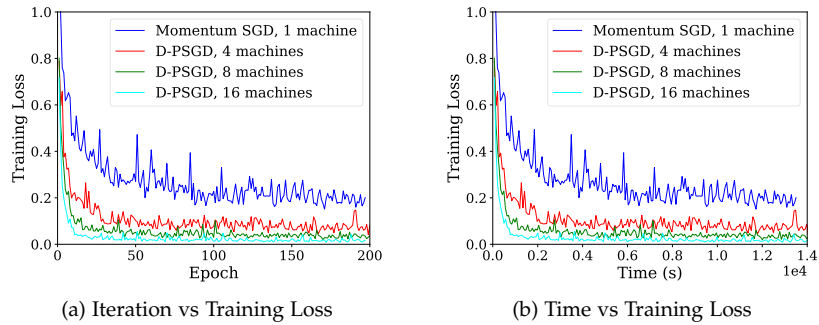

(a) Iteration vs Training Loss

(b) Time vs Training Loss

Figure 7: Training Loss comparison between D-PSGD on 1, 4, 8 and 16 machines

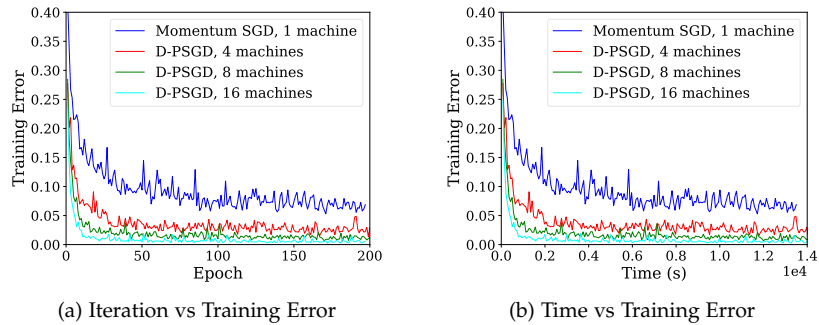

(a) Iteration vs Training Error

(b) Time vs Training Error

Figure 8: Training Error comparison between D-PSGD on 1, 4, 8 and 16 machines

Because our algorithm is synchronous, we desire each node to compute each minibatch roughly within the same time. If each machine has different computation power, we can use different minibatch sizes to compensate the speed difference, or allow faster machines to make more than 1 minibatch before synchronization.

# 6 Conclusion

This paper studies the D-PSGD algorithm on the decentralized computational network. We prove that D-PSGD achieves the same convergence rate (or equivalently computational complexity) as the C-PSGD algorithm, but outperforms C-PSGD by avoiding the communication traffic jam. To the best of our knowledge, this is the first work to show that decentralized algorithms admit the linear speedup and can outperform centralized algorithms.

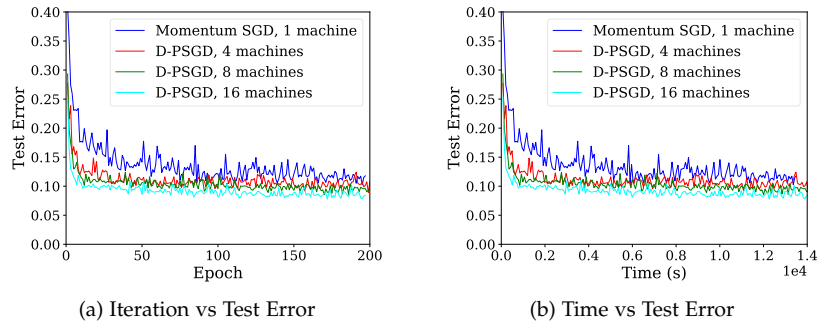

(a) Iteration vs Test Error

(b) Time vs Test Error

Figure 9: Test Error comparison between D-PSGD on 1, 4, 8 and 16 machines

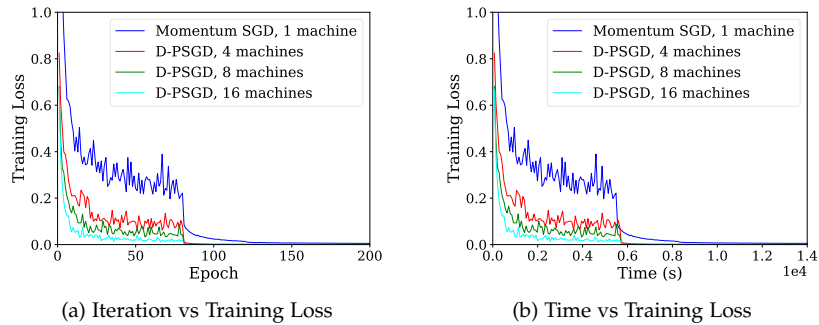

(a) Iteration vs Training Loss

(b) Time vs Training Loss

Figure 10: Training Loss comparison between D-PSGD on 1, 4, 8 and 16 machines

**Limitation and Future Work** The potential limitation of D-PSGD lies on the cost of synchronization. Breaking the synchronization barrier could make the decentralize algorithms even more efficient, but requires more complicated analysis. We will leave this direction for the future work.

On the system side, one future direction is to deploy D-PSGD to larger clusters beyond 112 GPUs and one such environment is state-of-the-art supercomputers. In such environment, we envision D-PSGD to be one necessary building blocks for multiple "centralized groups" to communicate. It is also interesting to deploy D-PSGD to mobile environments.

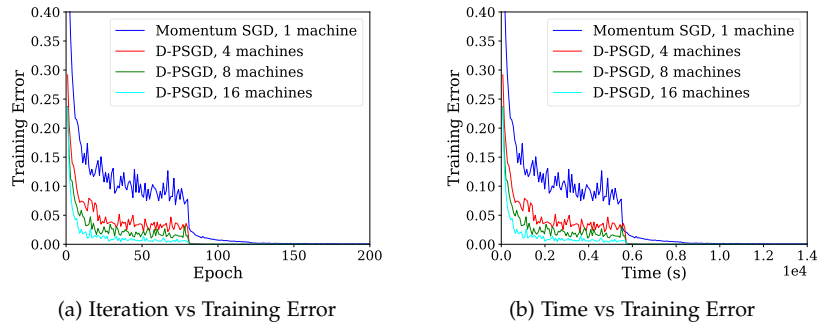

Figure 11: Training Error comparison between D-PSGD on 1, 4, 8 and 16 machines

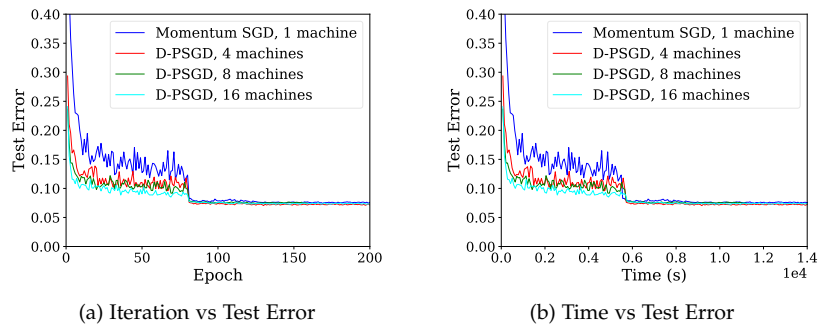

Figure 12: Test Error comparison between D-PSGD on 1, 4, 8 and 16 machines

# Supplemental Materials: More Experiments

## Industrial benchmark

In this section, we evaluate the effectiveness of our algorithm on IBM Watson Natural Language Classifier (NLC) workload. IBM Watson Natural Language Classifier (NLC) service, IBM's most popular cognitive service offering, is used by thousands of enterprise-level clients around the globe. The NLC task is to classify input sentences into a target category in a predefined label set. NLC has been extensively used in many practical applications, including sentiment analysis, topic classification, and question classification. At the core of NLC training is a CNN model that has a word-embedding lookup table layer, a convolutional layer and a fully connected layer with a softmax output layer. NLC is implemented using the Torch open-source deep learning framework.

**Methodology**  We use two datasets in our evaluation. The first dataset Joule is an in-house customer dataset that has 2.5K training samples, 1K test samples, and 311 different classes. The second dataset Yelp, which is a public dataset, has 500K training samples, 2K test samples and 5 different classes. The experiments are conducted on an IBM Power server, which has 40 IBM P8 cores, each core is 4-way SMP with clock frequence of 2GHz. The server has 128GB memory and is equipped with 8 K80 GPUs. DataParallelTable (DPT) is a NCCL-basedNvidia module in Torch that can leverage multiple GPUs to carry out centralized parallel SGD algorithm. NCCL is an all-reduce based implementation. We implemented the decentralized SGD algorithm in the NLC product. We now compare the convergence rate of centralized SGD (i.e. DPT) and our decentralized SGD implementation.

**Convergence results and test accuracy**  First, we examine the Joule dataset. We use 8 nodes and each node calculates with a mini-batch size of 2 and the entire run passes through 200 epochs. Figure 13 shows that centralized SGD algorithm and decentralized SGD algorithm achieve similar training loss (0.96) at roughly same convergence rate. Figure 14 shows that centralized SGD algorithm and decentralized SGD algorithm achieve similar testing error (43%). In the meantime, the communication cost is reduced by 3X in decentralized SGD case compared to the centralized SGD algorithm. Second, we examine the Yelp dataset. We use 8 nodes and each node calculates with a mini-batch size of 32 and the entire run passes through 20 epochs. Figure 13 shows that centralized SGD algorithm and decentralized SGD algorithm achieve similar training loss (0.86). Figure 14 shows that centralized SGD algorithm and decentralized SGD algorithm achieve similar testing error (39%). In the meantime, the communication cost is reduced by 2X in decentralized SGD case compared to the decentralized SGD case.

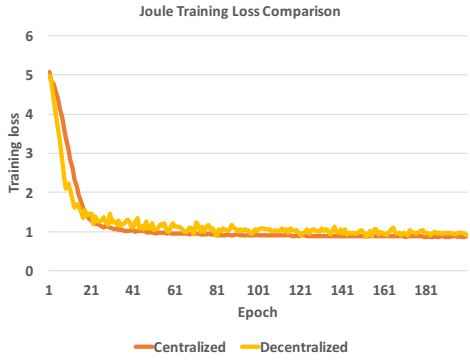

Figure 13: Training loss on Joule dataset

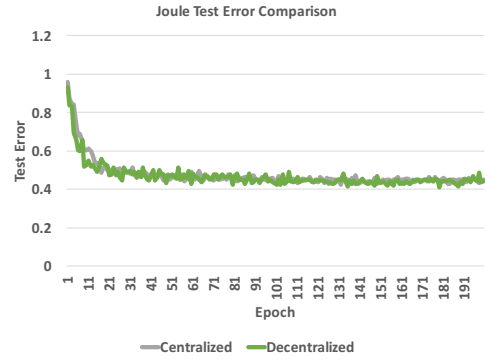

Figure 14: Test error on Joule dataset

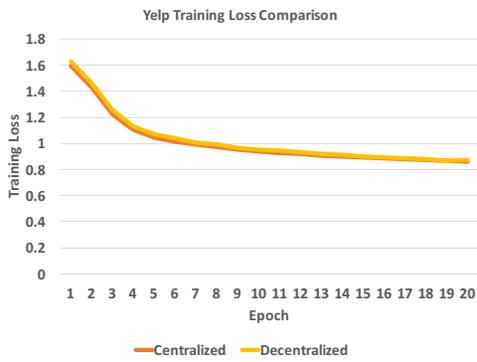

Figure 15: Training loss on Yelp dataset

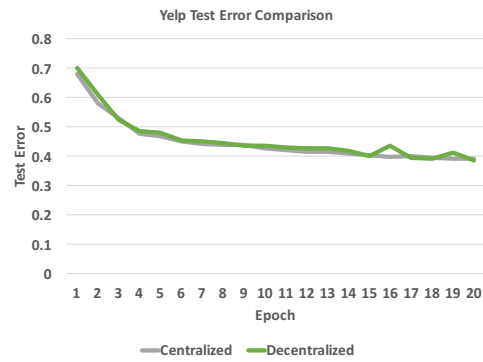

Figure 16: Test error on Yelp dataset

# Supplemental Materials: Proofs

We provide the proof to all theoretical results in this paper in this section.

**Lemma 5.** *Under Assumption 1 we have*

$$\left\|\frac{\mathbf{1}_n}{n} - W^k e_i\right\|^2 \leq \rho^k, \quad \forall i \in \{1, 2, \ldots, n\}, k \in \mathbb{N}.$$

*Proof.* Let $W^\infty := \lim_{k \to \infty} W^k$. Note that from Assumption 1-2 we have $\frac{\mathbf{1}_n}{n} = W^\infty e_i, \forall i$ since $W$ is doubly stochastic and $\rho < 1$. Thus

$$\begin{aligned}
\left\|\frac{\mathbf{1}_n}{n} - W^k e_i\right\|^2 &= \|(W^\infty - W^k)e_i\|^2 \\
&\leq \|W^\infty - W^k\|^2 \|e_i\|^2 \\
&= \|W^\infty - W^k\|^2 \\
&\leq \rho^k,
\end{aligned}$$

where the last step comes from the diagonalizability of $W$, completing the proof. □

**Lemma 6.** *We have the following inequality under Assumption 1:*

$$\mathbb{E}\|\partial f(X_j)\|^2 \leq \sum_{h=1}^{n} 3\mathbb{E}L^2 \left\|\frac{\sum_{i'=1}^{n} x_{j,i'}}{n} - x_{j,h}\right\|^2 + 3n\varsigma^2 + 3\mathbb{E}\left\|\nabla f\left(\frac{X_j \mathbf{1}_n}{n}\right)\mathbf{1}_n^\top\right\|^2, \forall j.$$

*Proof.* We consider the upper bound of $\mathbb{E}\|\partial f(X_j)\|^2$ in the following:

$$\begin{aligned}
&\mathbb{E}\|\partial f(X_j)\|^2 \\
&\leq 3\mathbb{E}\left\|\partial f(X_j) - \partial f\left(\frac{X_j \mathbf{1}_n}{n}\mathbf{1}_n^\top\right)\right\|^2 \\
&\quad + 3\mathbb{E}\left\|\partial f\left(\frac{X_j \mathbf{1}_n}{n}\mathbf{1}_n^\top\right) - \nabla f\left(\frac{X_j \mathbf{1}_n}{n}\right)\mathbf{1}_n^\top\right\|^2 \\
&\quad + 3\mathbb{E}\left\|\nabla f\left(\frac{X_j \mathbf{1}_n}{n}\right)\mathbf{1}_n^\top\right\|^2 \\
&\stackrel{\text{(Assumption 1-3)}}{\leq} 3\mathbb{E}\left\|\partial f(X_j) - \partial f\left(\frac{X_j \mathbf{1}_n}{n}\mathbf{1}_n^\top\right)\right\|_F^2 \\
&\quad + 3n\varsigma^2 \\
&\quad + 3\mathbb{E}\left\|\nabla f\left(\frac{X_j \mathbf{1}_n}{n}\right)\mathbf{1}_n^\top\right\|^2 \\
&\stackrel{\text{(Assumption 1-1)}}{\leq} \sum_{h=1}^{n} 3\mathbb{E}L^2 \left\|\frac{\sum_{i'=1}^{n} x_{j,i'}}{n} - x_{j,h}\right\|^2 + 3n\varsigma^2 + 3\mathbb{E}\left\|\nabla f\left(\frac{X_j \mathbf{1}_n}{n}\right)\mathbf{1}_n^\top\right\|^2.
\end{aligned}$$

This completes the proof. □

**Proof to Theorem 1.** We start form $f\left(\frac{X_{k+1}\mathbf{1}_n}{n}\right)$:

$$\mathbb{E}f\left(\frac{X_{k+1}\mathbf{1}_n}{n}\right)$$
$$=\mathbb{E}f\left(\frac{X_k W\mathbf{1}_n}{n}-\gamma\frac{\partial F(X_k;\xi_k)\mathbf{1}_n}{n}\right)$$
$$\stackrel{\text{(Assumption 1-2)}}{=}\mathbb{E}f\left(\frac{X_k\mathbf{1}_n}{n}-\gamma\frac{\partial F(X_k;\xi_k)\mathbf{1}_n}{n}\right)$$
$$\leqslant \mathbb{E}f\left(\frac{X_k\mathbf{1}_n}{n}\right)-\gamma\mathbb{E}\left\langle \nabla f\left(\frac{X_k\mathbf{1}_n}{n}\right),\frac{\partial f(X_k)\mathbf{1}_n}{n}\right\rangle$$
$$+\frac{\gamma^2 L}{2}\mathbb{E}\left\|\frac{\sum_{i=1}^n \nabla F_i(x_{k,i};\xi_{k,i})}{n}\right\|^2. \tag{7}$$

Note that for the last term we can split it into two terms:

$$\mathbb{E}\left\|\frac{\sum_{i=1}^n \nabla F_i(x_{k,i};\xi_{k,i})}{n}\right\|^2 = \mathbb{E}\left\|\frac{\sum_{i=1}^n \nabla F_i(x_{k,i};\xi_{k,i})-\sum_{i=1}^n \nabla f_i(x_{k,i})}{n}+\frac{\sum_{i=1}^n \nabla f_i(x_{k,i})}{n}\right\|^2$$
$$=\mathbb{E}\left\|\frac{\sum_{i=1}^n \nabla F_i(x_{k,i};\xi_{k,i})-\sum_{i=1}^n \nabla f_i(x_{k,i})}{n}\right\|^2$$
$$+\mathbb{E}\left\|\frac{\sum_{i=1}^n \nabla f_i(x_{k,i})}{n}\right\|^2$$
$$+\mathbb{E}\left\langle \frac{\sum_{i=1}^n \nabla F_i(x_{k,i};\xi_{k,i})-\sum_{i=1}^n \nabla f_i(x_{k,i})}{n},\frac{\sum_{i=1}^n \nabla f_i(x_{k,i})^2}{n}\right\rangle$$
$$=\mathbb{E}\left\|\frac{\sum_{i=1}^n \nabla F_i(x_{k,i};\xi_{k,i})-\sum_{i=1}^n \nabla f_i(x_{k,i})}{n}\right\|^2$$
$$+\mathbb{E}\left\|\frac{\sum_{i=1}^n \nabla f_i(x_{k,i})}{n}\right\|^2$$
$$+\mathbb{E}\left\langle \frac{\sum_{i=1}^n \mathbb{E}_{\xi_{k,i}}\nabla F_i(x_{k,i};\xi_{k,i})-\sum_{i=1}^n \nabla f_i(x_{k,i})}{n},\frac{\sum_{i=1}^n \nabla f_i(x_{k,i})^2}{n}\right\rangle$$
$$=\mathbb{E}\left\|\frac{\sum_{i=1}^n \nabla F_i(x_{k,i};\xi_{k,i})-\sum_{i=1}^n \nabla f_i(x_{k,i})}{n}\right\|^2$$
$$+\mathbb{E}\left\|\frac{\sum_{i=1}^n \nabla f_i(x_{k,i})}{n}\right\|^2.$$

Then it follows from (7) that

$$\mathbb{E}f\left(\frac{X_{k+1}\mathbf{1}_n}{n}\right)$$
$$\leqslant \mathbb{E}f\left(\frac{X_k\mathbf{1}_n}{n}\right)-\gamma\mathbb{E}\left\langle \nabla f\left(\frac{X_k\mathbf{1}_n}{n}\right),\frac{\partial f(X_k)\mathbf{1}_n}{n}\right\rangle$$
$$+\frac{\gamma^2 L}{2}\mathbb{E}\left\|\frac{\sum_{i=1}^n \nabla F_i(x_{k,i};\xi_{k,i})-\sum_{i=1}^n \nabla f_i(x_{k,i})}{n}\right\|^2$$

$$+ \frac{\gamma^2 L}{2} \mathbb{E} \left\| \frac{\sum_{i=1}^{n} \nabla f_i(x_{k,i})}{n} \right\|^2. \tag{8}$$

For the second last term we can bound it using $\sigma$:

$$\frac{\gamma^2 L}{2} \mathbb{E} \left\| \frac{\sum_{i=1}^{n} \nabla F_i(x_{k,i}; \xi_{k,i}) - \sum_{i=1}^{n} \nabla f_i(x_{k,i})}{n} \right\|^2$$

$$= \frac{\gamma^2 L}{2n^2} \sum_{i=1}^{n} \mathbb{E} \|\nabla F_i(x_{k,i}; \xi_{k,i}) - \nabla f_i(x_{k,i})\|^2$$

$$+ \frac{\gamma^2 L}{n^2} \sum_{i=1}^{n} \sum_{i'=i+1}^{n} \mathbb{E} \langle \nabla F_i(x_{k,i}; \xi_{k,i}) - \nabla f_i(x_{k,i}), \nabla F_{i'}(x_{k,i'}; \xi_{k,i'}) - \nabla f_{i'}(x_{k,i'}) \rangle$$

$$= \frac{\gamma^2 L}{2n^2} \sum_{i=1}^{n} \mathbb{E} \|\nabla F_i(x_{k,i}; \xi_{k,i}) - \nabla f_i(x_{k,i})\|^2$$

$$+ \frac{\gamma^2 L}{n^2} \sum_{i=1}^{n} \sum_{i'=i+1}^{n} \mathbb{E} \langle \nabla F_i(x_{k,i}; \xi_{k,i}) - \nabla f_i(x_{k,i}), \mathbb{E}_{\xi_{k,i'}} \nabla F_{i'}(x_{k,i'}; \xi_{k,i'}) - \nabla f_{i'}(x_{k,i'}) \rangle$$

$$= \frac{\gamma^2 L}{2n^2} \sum_{i=1}^{n} \mathbb{E} \|\nabla F_i(x_{k,i}; \xi_{k,i}) - \nabla f_i(x_{k,i})\|^2$$

$$\leqslant \frac{\gamma^2 L}{2n} \sigma^2,$$

where the last step comes from Assumption 1-3.

Thus it follows from (8):

$$\mathbb{E} f\left(\frac{X_{k+1} \mathbf{1}_n}{n}\right)$$

$$\leqslant \mathbb{E} f\left(\frac{X_k \mathbf{1}_n}{n}\right) - \gamma \mathbb{E} \left\langle \nabla f\left(\frac{X_k \mathbf{1}_n}{n}\right), \frac{\partial f(X_k) \mathbf{1}_n}{n} \right\rangle + \frac{\gamma^2 L}{2} \frac{\sigma^2}{n}$$

$$+ \frac{\gamma^2 L}{2} \mathbb{E} \left\| \frac{\sum_{i=1}^{n} \nabla f_i(x_{k,i})}{n} \right\|^2$$

$$= \mathbb{E} f\left(\frac{X_k \mathbf{1}_n}{n}\right) - \frac{\gamma - \gamma^2 L}{2} \mathbb{E} \left\| \frac{\partial f(X_k) \mathbf{1}_n}{n} \right\|^2 - \frac{\gamma}{2} \mathbb{E} \left\| \nabla f\left(\frac{X_k \mathbf{1}_n}{n}\right) \right\|^2 + \frac{\gamma^2 L}{2} \frac{\sigma^2}{n}$$

$$+ \frac{\gamma}{2} \underbrace{\mathbb{E} \left\| \nabla f\left(\frac{X_k \mathbf{1}_n}{n}\right) - \frac{\partial f(X_k) \mathbf{1}_n}{n} \right\|^2}_{=: T_1}, \tag{9}$$

where the last step comes from $2\langle a, b \rangle = \|a\|^2 + \|b\|^2 - \|a - b\|^2$.

We then bound $T_1$:

$$T_1 = \mathbb{E} \left\| \nabla f\left(\frac{X_k \mathbf{1}_n}{n}\right) - \frac{\partial f(X_k) \mathbf{1}_n}{n} \right\|^2$$

$$\leqslant \frac{1}{n} \sum_{i=1}^{n} \mathbb{E} \left\| \nabla f_i\left(\frac{\sum_{i'=1}^{n} x_{k,i'}}{n}\right) - \nabla f_i(x_{k,i}) \right\|^2$$

$$\overset{\text{(Assumption 1-1)}}{\leqslant} \frac{L^2}{n} \sum_{i=1}^{n} \underbrace{\mathbb{E} \left\| \frac{\sum_{i'=1}^{n} x_{k,i'}}{n} - x_{k,i} \right\|^2}_{=:Q_{k,i}}, \tag{10}$$

where we define $Q_{k,i}$ as the squared distance of the local optimization variable on the $i$-th node from the averaged local optimization variables on all nodes.

In order to bound $T_1$ we bound $Q_{k,i}$'s as the following:

$$\begin{aligned}
Q_{k,i} &= \mathbb{E} \left\| \frac{\sum_{i'=1}^{n} x_{k,i'}}{n} - x_{k,i} \right\|^2 \\
&= \mathbb{E} \left\| \frac{X_k \mathbf{1}_n}{n} - X_k e_i \right\|^2 \\
&= \mathbb{E} \left\| \frac{X_{k-1} W \mathbf{1}_n - \gamma \partial F(X_{k-1}; \xi_{k-1}) \mathbf{1}_n}{n} - (X_{k-1} W e_i - \gamma \partial F(X_{k-1}; \xi_{k-1}) e_i) \right\|^2 \\
&= \mathbb{E} \left\| \frac{X_{k-1} \mathbf{1}_n - \gamma \partial F(X_{k-1}; \xi_{k-1}) \mathbf{1}_n}{n} - (X_{k-1} W e_i - \gamma \partial F(X_{k-1}; \xi_{k-1}) e_i) \right\|^2 \\
&= \mathbb{E} \left\| \frac{X_0 \mathbf{1}_n - \sum_{i=0}^{k-1} \gamma \partial F(X_i; \xi_i) \mathbf{1}_n}{n} - \left( X_0 W^k e_i - \sum_{j=0}^{k-1} \gamma \partial F(X_j; \xi_j) W^{k-j-1} e_i \right) \right\|^2 \\
&= \mathbb{E} \left\| X_0 \left( \frac{\mathbf{1}_n}{n} - W^k e_i \right) - \sum_{j=0}^{k-1} \gamma \partial F(X_j; \xi_j) \left( \frac{\mathbf{1}_n}{n} - W^{k-j-1} e_i \right) \right\|^2 \\
&\overset{\text{(Assumption 1-4)}}{=} \mathbb{E} \left\| \sum_{j=0}^{k-1} \gamma \partial F(X_j; \xi_j) \left( \frac{\mathbf{1}_n}{n} - W^{k-j-1} e_i \right) \right\|^2 \\
&= \gamma^2 \mathbb{E} \left\| \sum_{j=0}^{k-1} \partial F(X_j; \xi_j) \left( \frac{\mathbf{1}_n}{n} - W^{k-j-1} e_i \right) \right\|^2 \\
&\leqslant 2\gamma^2 \underbrace{\mathbb{E} \left\| \sum_{j=0}^{k-1} (\partial F(X_j; \xi_j) - \partial f(X_j)) \left( \frac{\mathbf{1}_n}{n} - W^{k-j-1} e_i \right) \right\|^2}_{=:T_2} \\
&\quad + 2\gamma^2 \underbrace{\mathbb{E} \left\| \sum_{j=0}^{k-1} \partial f(X_j) \left( \frac{\mathbf{1}_n}{n} - W^{k-j-1} e_i \right) \right\|^2}_{=:T_3}.
\end{aligned} \tag{11}$$

For $T_2$, we provide the following upper bounds:

$$\begin{aligned}
T_2 &= \mathbb{E} \left\| \sum_{j=0}^{k-1} (\partial F(X_j; \xi_j) - \partial f(X_j)) \left( \frac{\mathbf{1}_n}{n} - W^{k-j-1} e_i \right) \right\|^2 \\
&= \sum_{j=0}^{k-1} \mathbb{E} \left\| (\partial F(X_j; \xi_j) - \partial f(X_j)) \left( \frac{\mathbf{1}_n}{n} - W^{k-j-1} e_i \right) \right\|^2
\end{aligned}$$

$$\leqslant \sum_{j=0}^{k-1} \mathbb{E}\|\partial F(X_j;\xi_j) - \partial f(X_j)\|^2 \left\|\frac{\mathbf{1}_n}{n} - W^{k-j-1}e_i\right\|^2$$

$$\leqslant \sum_{j=0}^{k-1} \mathbb{E}\|\partial F(X_j;\xi_j) - \partial f(X_j)\|_F^2 \left\|\frac{\mathbf{1}_n}{n} - W^{k-j-1}e_i\right\|^2$$

$$\overset{\text{(Lemma 5,Assumption 1-3)}}{\leqslant} n\sigma^2 \sum_{j=0}^{k-1} \rho^{k-j-1}$$

$$\leqslant \frac{n\sigma^2}{1-\rho}.$$

For $T_3$, we provide the following upper bounds:

$$T_3 = \mathbb{E}\left\|\sum_{j=0}^{k-1} \partial f(X_j)\left(\frac{\mathbf{1}_n}{n} - W^{k-j-1}e_i\right)\right\|^2$$

$$= \underbrace{\sum_{j=0}^{k-1} \mathbb{E}\left\|\partial f(X_j)\left(\frac{\mathbf{1}_n}{n} - W^{k-j-1}e_i\right)\right\|^2}_{=:T_4}$$

$$+ \underbrace{\sum_{j\neq j'} \mathbb{E}\left\langle \partial f(X_j)\left(\frac{\mathbf{1}_n}{n} - W^{k-j-1}e_i\right), \partial f(X_{j'})\left(\frac{\mathbf{1}_n}{n} - W^{k-j'-1}e_i\right)\right\rangle}_{=:T_5}$$

To bound $T_3$ we bound $T_4$ and $T_5$ in the following: for $T_4$,

$$T_4 = \sum_{j=0}^{k-1} \mathbb{E}\left\|\partial f(X_j)\left(\frac{\mathbf{1}_n}{n} - W^{k-j-1}e_i\right)\right\|^2$$

$$\leqslant \sum_{j=0}^{k-1} \mathbb{E}\|\partial f(X_j)\|^2 \left\|\frac{\mathbf{1}_n}{n} - W^{k-j}e_i\right\|^2$$

$$\overset{\text{(Lemmas 5 and 6)}}{\leqslant} 3\sum_{j=0}^{k-1}\sum_{h=1}^{n} \mathbb{E}L^2 Q_{j,h}\left\|\frac{\mathbf{1}_n}{n} - W^{k-j-1}e_i\right\|^2 + 3n\varsigma^2 \frac{1}{1-\rho}$$

$$+ 3\sum_{j=0}^{k-1} \mathbb{E}\left\|\nabla f\left(\frac{X_j\mathbf{1}_n}{n}\right)\mathbf{1}_n^\top\right\|^2 \left\|\frac{\mathbf{1}_n}{n} - W^{k-j-1}e_i\right\|^2.$$

We bound $T_5$ using two new terms $T_6$ and $T_7$:

$$T_5 = \sum_{j\neq j'}^{k-1} \mathbb{E}\left\langle \partial f(X_j)\left(\frac{\mathbf{1}_n}{n} - W^{k-j-1}e_i\right), \partial f(X_{j'})\left(\frac{\mathbf{1}_n}{n} - W^{k-j'-1}e_i\right)\right\rangle$$

$$\leqslant \sum_{j\neq j'}^{k-1} \mathbb{E}\left\|\partial f(X_j)\left(\frac{\mathbf{1}_n}{n} - W^{k-j-1}e_i\right)\right\| \left\|\partial f(X_{j'})\left(\frac{\mathbf{1}_n}{n} - W^{k-j'-1}e_i\right)\right\|$$

$$\leqslant \sum_{j\neq j'}^{k-1} \mathbb{E}\|\partial f(X_j)\| \left\|\frac{\mathbf{1}_n}{n} - W^{k-j-1}e_i\right\| \|\partial f(X_{j'})\| \left\|\frac{\mathbf{1}_n}{n} - W^{k-j'-1}e_i\right\|$$

$$\leqslant \sum_{j\neq j'}^{k-1} \mathbb{E}\|\partial f(X_j)\| \left\|\frac{\mathbf{1}_n}{n} - W^{k-j-1}e_i\right\| \|\partial f(X_{j'})\| \left\|\frac{\mathbf{1}_n}{n} - W^{k-j'-1}e_i\right\|$$

$$\leqslant \sum_{j\neq j'}^{k-1} \mathbb{E}\frac{\|\partial f(X_j)\|^2}{2}\left\|\frac{\mathbf{1}_n}{n}-W^{k-j-1}e_i\right\|\left\|\frac{\mathbf{1}_n}{n}-W^{k-j'-1}e_i\right\|$$

$$+\sum_{j\neq j'}^{k-1} \mathbb{E}\frac{\|\partial f(X_{j'})\|^2}{2}\left\|\frac{\mathbf{1}_n}{n}-W^{k-j-1}e_i\right\|\left\|\frac{\mathbf{1}_n}{n}-W^{k-j'-1}e_i\right\|$$

$$\stackrel{\text{Lemma 5}}{\leqslant} \sum_{j\neq j'}^{k-1} \mathbb{E}\left(\frac{\|\partial f(X_j)\|^2}{2}+\frac{\|\partial f(X_{j'})\|^2}{2}\right)\rho^{k-\frac{j+j'}{2}-1}$$

$$=\sum_{j\neq j'}^{k-1} \mathbb{E}(\|\partial f(X_j)\|^2)\rho^{k-\frac{j+j'}{2}-1}$$

$$\stackrel{\text{Lemma 6}}{\leqslant} \underbrace{3\sum_{j\neq j'}^{k-1}\left(\sum_{h=1}^{n}\mathbb{E}L^2 Q_{j,h}+\mathbb{E}\left\|\nabla f\left(\frac{X_j\mathbf{1}_n}{n}\right)\mathbf{1}_n^\top\right\|^2\right)\rho^{k-\frac{j+j'}{2}-1}}_{=:T_6}$$

$$+\underbrace{\sum_{j\neq j'}^{k-1} 3n\varsigma^2 \rho^{k-1-\frac{j+j'}{2}}}_{=:T_7},$$

where $T_7$ can be bounded using $\varsigma$ and $\rho$:

$$T_7 = 6n\varsigma^2 \sum_{j>j'}^{k-1}\rho^{k-1-\frac{j+j'}{2}}$$

$$=6n\varsigma^2 \frac{\left(\rho^{k/2}-1\right)\left(\rho^{k/2}-\sqrt{\rho}\right)}{(\sqrt{\rho}-1)^2(\sqrt{\rho}+1)}$$

$$\leqslant 6n\varsigma^2 \frac{1}{(1-\sqrt{\rho})^2},$$

and we bound $T_6$:

$$T_6 = 3\sum_{j\neq j'}^{k-1}\left(\sum_{h=1}^{n}\mathbb{E}L^2 Q_{j,h}+\mathbb{E}\left\|\nabla f\left(\frac{X_j\mathbf{1}_n}{n}\right)\mathbf{1}_n^\top\right\|^2\right)\rho^{k-\frac{j+j'}{2}-1}$$

$$=6\sum_{j=0}^{k-1}\left(\sum_{h=1}^{n}\mathbb{E}L^2 Q_{j,h}+\mathbb{E}\left\|\nabla f\left(\frac{X_j\mathbf{1}_n}{n}\right)\mathbf{1}_n^\top\right\|^2\right)\sum_{j'=j+1}^{k-1}\sqrt{\rho}^{2k-j-j'-2}$$

$$\leqslant 6\sum_{j=0}^{k-1}\left(\sum_{h=1}^{n}\mathbb{E}L^2 Q_{j,h}+\mathbb{E}\left\|\nabla f\left(\frac{X_j\mathbf{1}_n}{n}\right)\mathbf{1}_n^\top\right\|^2\right)\frac{\sqrt{\rho}^{k-j-1}}{1-\sqrt{\rho}}.$$

Plugging $T_6$ and $T_7$ into $T_5$ and then plugging $T_5$ and $T_4$ into $T_3$ yield the upper bound for $T_3$:

$$T_3 \leqslant 3\sum_{j=0}^{k-1}\sum_{h=1}^{n}\mathbb{E}L^2 Q_{j,h}\left\|\frac{\mathbf{1}_n}{n}-W^{k-j-1}e_i\right\|^2$$

$$+3\sum_{j=0}^{k-1}\mathbb{E}\left\|\nabla f\left(\frac{X_j\mathbf{1}_n}{n}\right)\mathbf{1}_n^\top\right\|^2\left\|\frac{\mathbf{1}_n}{n}-W^{k-j-1}e_i\right\|^2$$

$$
\begin{aligned}
&+ 6 \sum_{j=0}^{k-1} \left( \sum_{h=1}^{n} \mathbb{E} L^2 Q_{j,h} + \mathbb{E} \left\| \nabla f\left(\frac{X_j \mathbf{1}_n}{n}\right) \mathbf{1}_n^\top \right\|^2 \right) \frac{\sqrt{\rho}^{k-j-1}}{1-\sqrt{\rho}} \\
&+ \frac{3n\varsigma^2}{1-\rho} + \frac{6n\varsigma^2}{(1-\sqrt{\rho})^2} \\
\leq & 3 \sum_{j=0}^{k-1} \sum_{h=1}^{n} \mathbb{E} L^2 Q_{j,h} \left\| \frac{\mathbf{1}_n}{n} - W^{k-j-1} e_i \right\|^2 \\
&+ 3 \sum_{j=0}^{k-1} \mathbb{E} \left\| \nabla f\left(\frac{X_j \mathbf{1}_n}{n}\right) \mathbf{1}_n^\top \right\|^2 \left\| \frac{\mathbf{1}_n}{n} - W^{k-j-1} e_i \right\|^2 \\
&+ 6 \sum_{j=0}^{k-1} \left( \sum_{h=1}^{n} \mathbb{E} L^2 Q_{j,h} + \mathbb{E} \left\| \nabla f\left(\frac{X_j \mathbf{1}_n}{n}\right) \mathbf{1}_n^\top \right\|^2 \right) \frac{\sqrt{\rho}^{k-j-1}}{1-\sqrt{\rho}} \\
&+ \frac{9n\varsigma^2}{(1-\sqrt{\rho})^2},
\end{aligned}
$$

where the last step we use the fact that $\frac{1}{1-\rho} \leq \frac{1}{(1-\sqrt{\rho})^2}$.

Putting the bound for $T_2$ and $T_3$ back to (11) we get the bound for $Q_{k,i}$:

$$
\begin{aligned}
Q_{k,i} \leq & \frac{2\gamma^2 n \sigma^2}{1-\rho} + 6\gamma^2 \sum_{j=0}^{k-1} \sum_{h=1}^{n} \mathbb{E} L^2 \left\| \frac{\sum_{i'=1}^{n} x_{j,i'}}{n} - x_{j,h} \right\|^2 \left\| \frac{\mathbf{1}_n}{n} - W^{k-j-1} e_i \right\|^2 \\
&+ 6\gamma^2 \sum_{j=0}^{k-1} \mathbb{E} \left\| \nabla f\left(\frac{X_j \mathbf{1}_n}{n}\right) \mathbf{1}_n^\top \right\|^2 \left\| \frac{\mathbf{1}_n}{n} - W^{k-j-1} e_i \right\|^2 \\
&+ 12\gamma^2 \sum_{j=0}^{k-1} \left( \sum_{h=1}^{n} \mathbb{E} L^2 \left\| \frac{\sum_{i'=1}^{n} x_{j,i'}}{n} - x_{j,h} \right\|^2 + \mathbb{E} \left\| \nabla f\left(\frac{X_j \mathbf{1}_n}{n}\right) \mathbf{1}_n^\top \right\|^2 \right) \frac{\sqrt{\rho}^{k-j-1}}{1-\sqrt{\rho}} \\
&+ \frac{18\gamma^2 n \varsigma^2}{(1-\sqrt{\rho})^2} \\
\overset{\text{Lemma 5}}{\leq} & \frac{2\gamma^2 n \sigma^2}{1-\rho} + \frac{18\gamma^2 n \varsigma^2}{(1-\sqrt{\rho})^2} \\
&+ 6\gamma^2 \sum_{j=0}^{k-1} \sum_{h=1}^{n} \mathbb{E} L^2 Q_{j,h} \rho^{k-j-1} \\
&+ 6\gamma^2 \sum_{j=0}^{k-1} \mathbb{E} \left\| \nabla f\left(\frac{X_j \mathbf{1}_n}{n}\right) \mathbf{1}_n^\top \right\|^2 \rho^{k-j-1} \\
&+ 12\gamma^2 \sum_{j=0}^{k-1} \left( \sum_{h=1}^{n} \mathbb{E} L^2 Q_{j,h} + \mathbb{E} \left\| \nabla f\left(\frac{X_j \mathbf{1}_n}{n}\right) \mathbf{1}_n^\top \right\|^2 \right) \frac{\sqrt{\rho}^{k-j-1}}{1-\sqrt{\rho}} \\
= & \frac{2\gamma^2 n \sigma^2}{1-\rho} + \frac{18\gamma^2 n \varsigma^2}{(1-\sqrt{\rho})^2} \\
&+ 6\gamma^2 \sum_{j=0}^{k-1} \mathbb{E} \left\| \nabla f\left(\frac{X_j \mathbf{1}_n}{n}\right) \mathbf{1}_n^\top \right\|^2 \left( \rho^{k-j-1} + \frac{2\sqrt{\rho}^{k-j-1}}{1-\sqrt{\rho}} \right)
\end{aligned}
$$

$$+ 6\gamma^2 \sum_{j=0}^{k-1} \sum_{h=1}^{n} \mathbb{E} L^2 Q_{j,h} \left( \frac{2\sqrt{\rho}^{k-j-1}}{1 - \sqrt{\rho}} + \rho^{k-j-1} \right). \tag{12}$$

Till now, we have the bound for $Q_{k,i}$. We continue by bounding its average $M_k$ on all nodes, which is defined by:

$$\mathbb{E} M_k := \frac{\mathbb{E} \sum_{i=1}^{n} Q_{k,i}}{n} \tag{13}$$

$$\stackrel{(12)}{\leqslant} \frac{2\gamma^2 n\sigma^2}{1-\rho} + \frac{18\gamma^2 n\varsigma^2}{(1-\sqrt{\rho})^2}$$

$$+ 6\gamma^2 \sum_{j=0}^{k-1} \mathbb{E} \left\| \nabla f\left(\frac{X_j \mathbf{1}_n}{n}\right) \mathbf{1}_n^\top \right\|^2 \left( \rho^{k-j-1} + \frac{2\sqrt{\rho}^{k-j-1}}{1-\sqrt{\rho}} \right)$$

$$+ 6\gamma^2 n L^2 \sum_{j=0}^{k-1} \mathbb{E} M_j \left( \frac{2\sqrt{\rho}^{k-j-1}}{1-\sqrt{\rho}} + \rho^{k-j-1} \right).$$

Summing from $k = 0$ to $K - 1$ we get:

$$\sum_{k=0}^{K-1} \mathbb{E} M_k \leqslant \frac{2\gamma^2 n\sigma^2}{1-\rho} K + \frac{18\gamma^2 n\varsigma^2}{(1-\sqrt{\rho})^2} K$$

$$+ 6\gamma^2 \sum_{k=0}^{K-1} \sum_{j=0}^{k-1} \mathbb{E} \left\| \nabla f\left(\frac{X_j \mathbf{1}_n}{n}\right) \mathbf{1}_n^\top \right\|^2 \left( \rho^{k-j-1} + \frac{2\sqrt{\rho}^{k-j-1}}{1-\sqrt{\rho}} \right)$$

$$+ 6\gamma^2 n L^2 \sum_{k=0}^{K-1} \sum_{j=0}^{k-1} \mathbb{E} M_j \left( \frac{2\sqrt{\rho}^{k-j-1}}{1-\sqrt{\rho}} + \rho^{k-j-1} \right)$$

$$\leqslant \frac{2\gamma^2 n\sigma^2}{1-\rho} K + \frac{18\gamma^2 n\varsigma^2}{(1-\sqrt{\rho})^2} K$$

$$+ 6\gamma^2 \sum_{k=0}^{K-1} \mathbb{E} \left\| \nabla f\left(\frac{X_k \mathbf{1}_n}{n}\right) \mathbf{1}_n^\top \right\|^2 \left( \sum_{i=0}^{\infty} \rho^i + \frac{2 \sum_{i=0}^{\infty} \sqrt{\rho}^i}{1-\sqrt{\rho}} \right)$$

$$+ 6\gamma^2 n L^2 \sum_{k=0}^{K-1} \mathbb{E} M_k \left( \frac{2\sum_{i=0}^{\infty} \sqrt{\rho}^i}{1-\sqrt{\rho}} + \sum_{i=0}^{\infty} \rho^i \right)$$

$$\leqslant \frac{2\gamma^2 n\sigma^2}{1-\rho} K + \frac{18\gamma^2 n\varsigma^2}{(1-\sqrt{\rho})^2} K$$

$$+ \frac{18}{(1-\sqrt{\rho})^2} \gamma^2 \sum_{k=0}^{K-1} \mathbb{E} \left\| \nabla f\left(\frac{X_k \mathbf{1}_n}{n}\right) \mathbf{1}_n^\top \right\|^2$$

$$+ \frac{18}{(1-\sqrt{\rho})^2} \gamma^2 n L^2 \sum_{k=0}^{K-1} \mathbb{E} M_k,$$

where the second step comes from rearranging the summations and the last step comes from the summation of geometric sequences.

Simply by rearranging the terms we get the bound for the summation of $\mathbb{E} M_k$'s from $k = 0$ to $K - 1$:

$$\left(1 - \frac{18}{(1-\sqrt{\rho})^2}\gamma^2 nL^2\right) \sum_{k=0}^{K-1} \mathbb{E}M_k$$

$$\leqslant \frac{2\gamma^2 n\sigma^2}{1-\rho}K + \frac{18\gamma^2 n\varsigma^2}{(1-\sqrt{\rho})^2}K + \frac{18}{(1-\sqrt{\rho})}\gamma^2 \sum_{k=0}^{K-1} \mathbb{E}\left\|\nabla f\left(\frac{X_k \mathbf{1}_n}{n}\right)\mathbf{1}_n^\top\right\|^2$$

$$\implies \sum_{k=0}^{K-1} \mathbb{E}M_k \leqslant \frac{2\gamma^2 n\sigma^2}{(1-\rho)\left(1 - \frac{18}{(1-\sqrt{\rho})^2}\gamma^2 nL^2\right)}K + \frac{18\gamma^2 n\varsigma^2}{(1-\sqrt{\rho})^2\left(1 - \frac{18}{(1-\sqrt{\rho})^2}\gamma^2 nL^2\right)}K$$

$$+ \frac{18\gamma^2}{(1-\sqrt{\rho})^2\left(1 - \frac{18}{(1-\sqrt{\rho})^2}\gamma^2 nL^2\right)} \sum_{k=0}^{K-1} \mathbb{E}\left\|\nabla f\left(\frac{X_k \mathbf{1}_n}{n}\right)\mathbf{1}_n^\top\right\|^2. \qquad (14)$$

Recall (10) that $T_1$ can be bounded using $M_k$:

$$\mathbb{E}T_1 \leqslant \frac{L^2}{n}\sum_{i=1}^n \mathbb{E}Q_{k,i} = L^2 \mathbb{E}M_k. \qquad (15)$$

We are finally able to bound the error by combining all above. Starting from (9):

$$\mathbb{E}f\left(\frac{X_{k+1}\mathbf{1}_n}{n}\right) \leqslant \mathbb{E}f\left(\frac{X_k \mathbf{1}_n}{n}\right) - \frac{\gamma - \gamma^2 L}{2}\mathbb{E}\left\|\frac{\partial f(X_k)\mathbf{1}_n}{n}\right\|^2 - \frac{\gamma}{2}\mathbb{E}\left\|\nabla f\left(\frac{X_k \mathbf{1}_n}{n}\right)\right\|^2$$

$$+ \frac{\gamma^2 L}{2n}\sigma^2 + \frac{\gamma}{2}\mathbb{E}T_1$$

$$\stackrel{(15)}{\leqslant} \mathbb{E}f\left(\frac{X_k \mathbf{1}_n}{n}\right) - \frac{\gamma - \gamma^2 L}{2}\mathbb{E}\left\|\frac{\partial f(X_k)\mathbf{1}_n}{n}\right\|^2 - \frac{\gamma}{2}\mathbb{E}\left\|\nabla f\left(\frac{X_k \mathbf{1}_n}{n}\right)\right\|^2$$

$$+ \frac{\gamma^2 L}{2n}\sigma^2 + \frac{\gamma}{2}L^2 \mathbb{E}M_k.$$

Summing from $k=0$ to $k=K-1$ we get:

$$\frac{\gamma - \gamma^2 L}{2}\sum_{k=0}^{K-1}\mathbb{E}\left\|\frac{\partial f(X_k)\mathbf{1}_n}{n}\right\|^2 + \frac{\gamma}{2}\sum_{k=0}^{K-1}\mathbb{E}\left\|\nabla f\left(\frac{X_k \mathbf{1}_n}{n}\right)\right\|^2$$

$$\leqslant f(0) - f^* + \frac{\gamma^2 KL}{2n}\sigma^2 + \frac{\gamma}{2}L^2 \sum_{k=0}^{K-1}\mathbb{E}M_k$$

$$\stackrel{(14)}{\leqslant} f(0) - f^* + \frac{\gamma^2 KL}{2n}\sigma^2$$

$$+ \frac{\gamma}{2}L^2 \frac{2\gamma^2 n\sigma^2}{(1-\rho)\left(1 - \frac{18}{(1-\sqrt{\rho})^2}\gamma^2 nL^2\right)}K + \frac{\gamma}{2}L^2 \frac{18\gamma^2 n\varsigma^2}{(1-\sqrt{\rho})^2\left(1 - \frac{18}{(1-\sqrt{\rho})^2}\gamma^2 nL^2\right)}K$$

$$+ \frac{\gamma}{2}L^2 \frac{18\gamma^2}{(1-\sqrt{\rho})^2\left(1 - \frac{18}{(1-\sqrt{\rho})^2}\gamma^2 nL^2\right)} \sum_{k=0}^{K-1}\mathbb{E}\left\|\nabla f\left(\frac{X_k \mathbf{1}_n}{n}\right)\mathbf{1}_n^\top\right\|^2$$

$$
\begin{aligned}
= & f(0) - f^* + \frac{\gamma^2 K L}{2n}\sigma^2 \\
& + \frac{\gamma^3 L^2 n \sigma^2}{(1-\rho)\left(1 - \frac{18}{(1-\sqrt{\rho})^2}\gamma^2 n L^2\right)} K + \frac{9\gamma^3 L^2 n \varsigma^2}{(1-\sqrt{\rho})^2\left(1 - \frac{18}{(1-\sqrt{\rho})^2}\gamma^2 n L^2\right)} K \\
& + \frac{9n\gamma^3 L^2}{(1-\sqrt{\rho})^2\left(1 - \frac{18}{(1-\sqrt{\rho})^2}\gamma^2 n L^2\right)} \sum_{k=0}^{K-1} \mathbb{E}\left\|\nabla f\left(\frac{X_k \mathbf{1}_n}{n}\right)\right\|^2
\end{aligned}
$$

By rearranging the inequality above, we obtain:

$$
\begin{aligned}
\Longrightarrow & \frac{\frac{\gamma - \gamma^2 L}{2}\sum_{k=0}^{K-1}\mathbb{E}\left\|\frac{\partial f(X_k)\mathbf{1}_n}{n}\right\|^2 + \left(\frac{\gamma}{2} - \frac{9n\gamma^3 L^2}{(1-\sqrt{\rho})^2\left(1-\frac{18}{(1-\sqrt{\rho})^2}\gamma^2 n L^2\right)}\right)\sum_{k=0}^{K-1}\mathbb{E}\left\|\nabla f\left(\frac{X_k\mathbf{1}_n}{n}\right)\right\|^2}{\gamma K} \\
\leqslant & \frac{f(0) - f^*}{\gamma K} + \frac{\gamma L}{2n}\sigma^2 + \frac{\gamma^2 L^2 n \sigma^2}{(1-\rho)\left(1-\frac{18}{(1-\sqrt{\rho})^2}\gamma^2 n L^2\right)} + \frac{9\gamma^2 L^2 n \varsigma^2}{(1-\sqrt{\rho})^2\left(1-\frac{18}{(1-\sqrt{\rho})^2}\gamma^2 n L^2\right)}.
\end{aligned}
$$

which completes the proof. $\square$

*Proof to Corollary 2.* Substitute $\gamma = \frac{1}{2L + \sigma\sqrt{K/n}}$ into Theorem 1 and remove the $\left\|\frac{\partial f(X_k)\mathbf{1}_n}{n}\right\|^2$ terms on the LHS. We get

$$
\begin{aligned}
& \frac{D_1 \sum_{k=0}^{K-1}\mathbb{E}\left\|\nabla f\left(\frac{X_k\mathbf{1}_n}{n}\right)\right\|^2}{K} \\
\leqslant & \frac{2(f(0)-f^*)L}{K} + \frac{(f(0)-f^*)\sigma}{\sqrt{Kn}} + \frac{L\sigma^2}{4nL + 2\sigma\sqrt{Kn}} \\
& + \frac{L^2 n}{(2L+\sigma\sqrt{K/n})^2 D_2}\left(\frac{\sigma^2}{1-\rho} + \frac{9\varsigma^2}{(1-\sqrt{\rho})^2}\right) \\
\leqslant & \frac{2(f(0)-f^*)L}{K} + \frac{(f(0)-f^* + L/2)\sigma}{\sqrt{Kn}} \\
& + \frac{L^2 n}{(\sigma\sqrt{K/n})^2 D_2}\left(\frac{\sigma^2}{1-\rho} + \frac{9\varsigma^2}{(1-\sqrt{\rho})^2}\right). \quad (16)
\end{aligned}
$$

We first show $D_1$ and $D_2$ are approximately constants when (6) is satisfied.

$$
D_1 := \left(\frac{1}{2} - \frac{9\gamma^2 L^2 n}{(1-\sqrt{\rho})^2 D_2}\right), \quad D_2 := \left(1 - \frac{18\gamma^2}{(1-\sqrt{\rho})^2} n L^2\right).
$$

Note that

$$
\begin{aligned}
\gamma^2 \leqslant \frac{(1-\sqrt{\rho})^2}{36nL^2} & \Longrightarrow D_2 \geqslant 1/2, \\
\gamma^2 \leqslant \frac{(1-\sqrt{\rho})^2}{72 L^2 n} & \Longrightarrow D_1 \geqslant 1/4.
\end{aligned}
$$

Since
$$\gamma^2 \leqslant \frac{n}{\sigma^2 K},$$
as long as we have
$$\frac{n}{\sigma^2 K} \leqslant \frac{(1-\sqrt{\rho})^2}{36nL^2}$$
$$\frac{n}{\sigma^2 K} \leqslant \frac{(1-\sqrt{\rho})^2}{72L^2 n},$$
$D_2 \geqslant 1/2$ and $D_1 \geqslant 1/4$ will be satisfied. Solving above inequalities we get (6).

Now with (6) we can safely replace $D_1$ and $D_2$ in (17) with $1/4$ and $1/2$ respectively. Thus

$$\begin{aligned}
\frac{\sum_{k=0}^{K-1} \mathbb{E}\left\|\nabla f\left(\frac{X_k \mathbf{1}_n}{n}\right)\right\|^2}{4K} \\
\leqslant \frac{2(f(0)-f^*)L}{K} + \frac{(f(0)-f^*+L/2)\sigma}{\sqrt{Kn}} \\
+ \frac{2L^2 n}{(\sigma\sqrt{K/n})^2}\left(\frac{\sigma^2}{1-\rho} + \frac{9\varsigma^2}{(1-\sqrt{\rho})^2}\right).
\end{aligned} \tag{17}$$

Given (5), the last term is bounded by the second term, completing the proof. □

*Proof to Theorem 3.* This can be seen from a simple analysis that the $\rho, \sqrt{\rho}$ for this $W$ are asymptotically $1 - \frac{16\pi^2}{3n^2}, 1 - \frac{8\pi^2}{3n^2}$ respectively when $n$ is large. Then by requiring (6) we need $n \leq O(K^{1/6})$. To satisfy (5) we need $n \leq O\left(K^{1/9}\right)$ when $\varsigma = 0$ and $n \leq O(K^{1/13})$ when $\varsigma > 0$. This completes the proof. □

*Proof to Theorem 4.* From (14) with $\gamma = \frac{1}{2L + \sigma\sqrt{K/n}}$ we have

$$\begin{aligned}
\frac{\sum_{k=0}^{K-1} \mathbb{E} M_k}{K} &\leqslant \frac{2\gamma^2 n \sigma^2}{(1-\rho)D_2} + \frac{18\gamma^2 n \varsigma^2}{(1-\sqrt{\rho})^2 D_2} \\
&\quad + \frac{18\gamma^2}{(1-\sqrt{\rho})^2 D_2} \frac{\sum_{k=0}^{K-1}\mathbb{E}\left\|\nabla f\left(\frac{X_k \mathbf{1}_n}{n}\right)\mathbf{1}_n^\top\right\|^2}{K} \\
&= \frac{2\gamma^2 n \sigma^2}{(1-\rho)D_2} + \frac{18\gamma^2 n \varsigma^2}{(1-\sqrt{\rho})^2 D_2} \\
&\quad + \frac{18\gamma^2 n}{(1-\sqrt{\rho})^2 D_2}\frac{\sum_{k=0}^{K-1}\mathbb{E}\left\|\nabla f\left(\frac{X_k \mathbf{1}_n}{n}\right)\right\|^2}{K} \\
&\overset{\text{Corollary 2}}{\leqslant} \frac{2\gamma^2 n\sigma^2}{(1-\rho)D_2} + \frac{18\gamma^2 n\varsigma^2}{(1-\sqrt{\rho})^2 D_2} + \frac{\gamma^2 L^2 n}{D_1 D_2}\left(\frac{\sigma^2}{1-\rho} + \frac{9\varsigma^2}{(1-\sqrt{\rho})^2}\right) \\
&\quad + \frac{18\gamma^2 n}{(1-\sqrt{\rho})^2 D_2}\left(\frac{f(0)-f^*}{\gamma K} + \frac{\gamma L \sigma^2}{2nD_1}\right)
\end{aligned}$$

$$= \frac{n\gamma^2}{D_2} A.$$

This completes the proof. □